\documentclass[12pt]{amsart}

\newtheorem{theorem}{Theorem}[section]

\newtheorem{corollary}[theorem]{Corollary}

\newtheorem{lemma}[theorem]{Lemma}

\theoremstyle{definition}

\newtheorem{remark}[theorem]{Remark}

\begin{document}

\title[]{On high dimensional maximal operators}

\author{J. M. Aldaz and J. P\'erez L\'azaro}
\address{Departamento de Matem\'aticas,
Universidad  Aut\'onoma de Madrid, Cantoblanco 28049, Madrid, Spain.}
\email{jesus.munarriz@uam.es}
\address{Departamento de Matem\'aticas y Computaci\'on,
Universidad  de La Rioja, 26004 Logro\~no, La Rioja, Spain.}
\email{javier.perezl@unirioja.es}


\thanks{The authors were partially supported by Grants MTM2012-37436-C02-02 and MTM2012-36732-C03-02 of the
D.G.I. of Spain, respectively.}







\begin{abstract} In this note we describe some recent advances in the area of maximal
function inequalities. We also study the behaviour of
 the centered Hardy-Littlewood maximal operator
associated to certain families of doubling, radial decreasing  measures,  
and acting on radial functions. In fact, we precisely determine when 
the weak type $(1,1)$ bounds are  uniform in the dimension.
\end{abstract}


\maketitle


\section {Introduction}

\markboth{J. M. Aldaz and J. P\'erez-L\'azaro}{High dimensional maximal operators}

Given a Borel measure $\mu$ on a metric space $X$  and a locally integrable
function $g$, the centered Hardy-Littlewood maximal operator $M_{\mu}$ is given by
\begin{equation}\label{HLMF}
M_{\mu} g(x) := \sup _{\{r > 0: 0 < \mu (B(x, r)) \}} \frac{1}{\mu
(B(x, r))} \int _{B(x, r)} \vert g\vert d\mu,
\end{equation}
where $B(x, r)$ denotes the  open ball of radius
$r > 0$ centered at $x$. Recall that $g$ is locally integrable if for every $x\in X$ there exists
an $r > 0$ such that $\int _{B(x, r)} \vert g\vert d\mu < \infty$. For instance, $g(x) := 1/x$ is locally
integrable on $(0,\infty)$, but not on $\mathbb{R}$, regardless of how it is extended to
$(-\infty, 0]$.

We allow measures that assign infinite size to some balls. Of course, if $\mu$ 
 assigns infinite measure
to all balls, then it is  of no interest in this context, since then $M_{\mu} g\equiv 0 $ for every locally integrable $g$
(we adopt the convention $\infty/\infty = \infty \cdot 0 = 0$).
Note that if all balls (with finite radii) have finite measure, then it does not matter 
whether one uses open or closed balls in the definition of $M_\mu$. It follows from 
countable additivity that this
does not alter the value of $M_{\mu} g(x)$, since closed  (resp. open) balls can be obtained as 
countable intersections (resp. unions) 
of open (resp. closed) balls with the same center.
   When  $\mu = \lambda^d$, the  $d$-dimensional Lebesgue measure,
we often simplify notation, by writing $M$ rather than $M_{\lambda^d}$
and $dx$ instead of $d \lambda^d (x)$.

It is well known that $M_\mu$ is a positive, sublinear operator, acting  on
the cone of positive, locally integrable functions
($M_\mu$ is defined by using $|g|$ rather than $g$). The Hardy-Littlewood maximal operator admits many variants: Instead of
 averaging $|g|$ over balls centered at
$x$ (the centered operator) as in (\ref{HLMF}), it is possible to consider all balls containing $x$ (the uncentered
operator)  or average over  convex bodies more general than
euclidean balls (and even over more
general sets, for instance, star-shaped, lower dimensional, etc.). It can also be applied to locally finite
measures $\nu$ (rather than just functions) by setting (say, in the centered case)
\begin{equation}\label{HLMFmeas}
M_{\mu} \nu(x) := \sup _{\{r > 0: \mu (B(x, r)) > 0\}} \frac{\nu (B(x, r))}{\mu
(B(x, r))} .
\end{equation}

The Hardy-Littlewood maximal operator is  an often used tool
 in  Real and Harmonic Analysis,
mainly (but not exclusively) due to the  fact that while
$|g|\le M_\mu g$ a.e.,
  $M_\mu g$  is
 not too large (in an $L^p$ sense)  since for every Borel measure $\mu$ defined on $\mathbb{R}^d$,
  it satisfies the following strong type $(p,p)$ inequality: 
$\|M_\mu g\|_p \le C_p \|g\|_p$ for  $1 < p \le \infty$.  
Thus,  $M_\mu g$ is often used
 to replace $g$, or some average of $g$, in chains of inequalities, without leaving $L^p$ ($p >1$).

 The situation when $p=1$ is different. Taking $g=\chi_{[0,1]}$, we see that $Mg$ (on the real line with Lebesgue
 measure) behaves essentially like $1/x$ near infinity, so $Mg$ is not integrable. However, 
 it follows from the Besicovitch Covering Theorem that
  $M_\mu$  satisfies the weak type $(1,1)$ inequality
 $\sup_{\alpha > 0}\alpha \mu (\{M_\mu g \ge \alpha\}) \le c_1 \|g\|_1$  for every Borel measure 
 $\mu$  on $\mathbb{R}^d$.
  This is a very important fact, as
 it implies the $L^p$ bounds for $1 < p < \infty$ via interpolation (the Marcinkiewicz Interpolation
Theorem generalizes this result).
 From now on we shall use  $c_{1,d}$  to denote the lowest possible constant in the weak type (1,1) inequality
 when the dimension is $d$,
and likewise,   $C_{p,d}$ will denote the lowest strong $(p,p)$ constant in 
dimension $d$.
 
 \section{Weak bounds, strong bounds, and dimensions}
 
 An aspect of the Hardy-Littlewood maximal operator that is receiving increasing
 attention, but which will not be touched upon here, is that of its regularity properties (cf. for instance
 \cite{AlPe1},  \cite{AlPe2},  \cite{AlPe3}, \cite{ACP} and the references contained
 therein). In this paper we restrict our attention to results regarding weak and strong type bounds.
 Since, as mentioned above, maximal operators are often used in chains of inequalities, improvements in
 these bounds lead to improvements in several other inequalities.

Considerable efforts
have gone into determining how changing the dimension of $\mathbb{R}^d$ modifies the
best constants $C_{p,d}$ and $c_{1,d}$ in the case of Lebesgue measure.
When $p =\infty$, we can take $C_{p,d} = 1$ in every dimension, since
averages never exceed a supremum. At the other endpoint $p=1$, the first boundedness
arguments used the Vitali covering lemma, which leads to exponential bounds of the type $c_{1,d} \le 3^d$,
and by interpolation, to exponential bounds for  $C_{p,d}$. So it is natural to try to improve on these bounds,
and in particular, to seek bounds independent of the dimension, with a view towards infinite dimensional
generalizations of Harmonic Analysis. 

In the Vitali covering lemma one obtains a disjoint subfamily from
a finite family of balls by a greedy algorithm and enlarging radii: Choose first the ball $B_1$ with largest
radius. Then remove from the collection all the balls that intersect it. Observe that the union of
these balls is contained in the ball $3 B_1$ with the same center and three times the radius as $B_1$.
Then choose $B_2$ as the ball with the largest radius among the balls left, and repeat. This argument
works well whenever the measure of balls with large radii is controlled by the measure of balls
with the same center and smaller radii, in the following sense: There exists  a constant $K$ such
that for all balls $B$, $\mu 2B \le K \mu B$. Such measures $\mu$ are called doubling because we
double the radius of $B$, but in fact any other constant $t >1$ could be used in place of 2. For
instance, doubling with 2 implies doubling with 4, with constant $K^2$, and  doubling
with 4 implies doubling with 2, trivially.

In his Princeton Ph. D. thesis, motivated by Fritz John's solution of the wave equation via spherical means, 
Prof. Antonio Cordoba 
(personal communication)
considered what nowadays is called Bourgain's
circular maximal function, where averages are taken over circumferences centered at a
point, in dimension $d=2$ (there is a
small subtlety in the definition; since circumferences have area zero, one needs to work
first with functions defined everywhere, for instance, continuous functions, or $C^\infty$ functions, and then, if
one manages to prove strong type bounds of some sort, the operator can be defined
over measurable functions via approximation arguments).  However, A. Cordoba was unable 
to obtain $L^p$ bounds for this maximal operator. As it turns out, these bounds were easier to
establish in higher dimensions. E. M. Stein showed that for $d \ge 3$ the (Stein's) spherical
maximal operator (where averages are taken over centered spheres) was bounded in $L^p$ if and 
only if $p > d/(d-1)$, cf. \cite{StPNAS}. It took about ten years, and the efforts of J. Bourgain,
to extend Stein's result to $d = 2$, cf. \cite{Bou0}. So the moral here seems to be that one should not
start with the hardest case. Of course, a priori it may not be obvious what is easy and what is
difficult. For instance, in $d = 1$  a simple covering argument yields, for the uncentered operator
and essentially all 
 measures, $c_{1,1} \le 2$ (cf. Theorem \ref{1d} below)
and often $c_{1,1} = 2$ is sharp (example: Lebesgue measure). 
However, if we ask the same question for the centered operator and (just) Lebesgue measure,
then even proving that the constant is different from 2 is difficult. This was done in
\cite{A1}, where the then commonly accepted conjecture $c_{1,1} = 3/2$ was also refuted.
The exact value $c_{1,1} = (11 + \sqrt{61})/12$ was obtained by A. Melas by a rather
involved argument,
 in the two papers \cite{Me1}, \cite{Me2}.

Returning to the spherical maximal operator, it is more or less intuitively clear that it controls the
Hardy-Littlewood maximal operator $M$ associated to euclidean balls (but this requires some argument). 
By proving dimension
independent bounds for the spherical maximal operator,
E. M. Stein showed that for $M$,
there exist bounds for $C_p$ that are independent of $d$ (\cite{St1}, \cite{St2},
\cite{StSt}, see also \cite{St3}).
Stein's result was
generalized to the maximal function defined using an arbitrary norm by J. Bourgain (\cite{Bou1},
\cite{Bou2}, \cite{Bou3}) and A. Carbery (\cite{Ca}) when $p>3/2$.
For $\ell_q$ balls, $1\le q <\infty$, D. M\"{u}ller \cite{Mu}
showed that uniform bounds again hold for every $p
> 1$ (given $1\le q <\infty$, the $\ell_q$ balls
are defined using the norm $\|x\|_q :=\left( |x_1|^q+ |x_2|^q+\dots + |x_d|^q\right)^{1/q}$). 

Regarding weak type $(1,1)$ inequalities, in \cite{StSt} E. M. Stein
and J. O. Str\"{o}mberg proved that the smallest
constants in the weak type (1,1) inequality satisfied by $M$ grow
at most like $O(d)$ for euclidean balls, using the heat semigroup, and at most like $O( d\log d)$
for more general balls, by a difficult covering lemma argument. 
They also asked if uniform bounds could be found,
a {\em question  still open} for euclidean balls. 

Semigroup theory enters maximal function estimates via the Hopf maximal ergodic
theorem for semigroups of operators, applied to the heat semigroup. 
Here the supremum is taken over time (one dimensional)
so the bound is independent of dimension. Now the maximal function bound $C d$ ($C$ a constant) appears
as follows:
It is possible to express the centered maximal operator in terms of convolutions:
$$
Mf(x) = \sup_{r > 0} |f|*\frac{\chi_{B(0,r)}}{\lambda^d(B(0,r))}(x).
$$
The argument then proceeds by showing that there exists a constant $C > 0$
and $s= s(d)$ such that
$$ 
\frac{\chi_{B(0,r)}}{\lambda^d(B(0,r))}(x) \le \frac{C d}{s}\int_0^s \frac{1}{(4\pi t)^{d/2}} e^{-\frac{\|x\|_2^2}{4t}} dt.
$$

These results about the Hardy-Littlewood maximal operator were obtained during the eighties, 
after which activity in this area slowed down. But recently, it seems to have picked up steam.
In 2008 the note \cite{A2} was posted in the Math ArXiv (but was published in 2011, so some 
papers that cite it have earlier publication dates). 
It is shown there that if one considers
cubes with sides parallel to the coordinate axes (that is,  $\ell_\infty$ balls)
 instead of euclidean balls, then the best constants $c_{1,d}$ must diverge to infinity with $d$, 
and thus the  answer
to the Stein-Str\"{o}mberg question
is negative for cubes.  This was proven by elementary means, basically calculus and first year probability 
(the normal approximation to the
binomial distribution). More advanced probabilistic techniques (the theory of stochastic processes and
in particular, the brownian bridge) quickly lead to an improvement:
 G. Aubrun 
 showed shortly after that $c_{1,d} \ge \Theta(\log^{1 -\varepsilon}d)$,
where $\Theta$ denotes the exact order and $\varepsilon > 0$ is arbitrary, cf. \cite{Au}. 
Finally,  
the question
whether
 the
maximal operator associated to cubes and Lebesgue measure is uniformly bounded in $d$,
for each
$1 < p \le 3/2$, has recently received a  positive 
answer  by J. Bourgain  (Math. ArXiv, December 11th, 2012). So, save for
refinements on the size of the constants, the situation is now well understood
for cubes (and Lebesgue measure).

These results  suggest (at least to us) that uniform bounds for $c_{1,d}$ may fail to exist if one uses euclidean balls (the original question of Stein and Str\"ombeg) since there seems to be
 no reason to believe that the maximal operator
associated to euclidean balls is substantially smaller than the maximal operator
associated to cubes.

A very significant extension of the Stein and Str\"{o}mberg's
$O( d\log d)$ theorem, beyond $\mathbb{R}^d$, has recently been  obtained  by A. Naor and T. Tao, cf. \cite{NaTa}. At 
the level of generality these authors work, the order of growth $O( d\log d)$ cannot be lowered, 
as they show by constructing
the appropriate counterexample.

 In the Vitali covering lemma one covers balls by expanding the radius
of an intersecting ball, which may have only slightly larger radius than the others. It was already noted in 
\cite{StSt} that engulfing balls by expanding the radius of a much larger ball 
can be more efficient. This idea leads Naor and Tao to define the Microdoubling and Strong Microdoubling
properties on metric measure spaces.

A metric measure space $(X, d,\mu)$ is a separable metric space $(X, d)$, equipped with a Radon
measure $\mu$. Naor and Tao also assume that $0 < \mu(B(x, r)) < \infty$ for all
$r > 0$.
Now  $(X, d,\mu)$ is definend to be $d$-Microdoubling with constant $K$ if for all $x \in X$ and all $r >0$,
we have
$$
\mu B\left(x,\left(1 + \frac{1}{d} \right)r\right) \le K \mu B(x,r).
$$ 
Note that the case $n = 1$ is just doubling.
And $(X, d,\mu)$ is  Strong $d$-Microdoubling  with constant $K$ if
for all $x$, all $r > 0$ and all $y\in B(x,r)$, 
$$
\mu B\left(y,\left(1 + \frac{1}{d} \right)r\right) \le K \mu B(x,r).
$$ 
Naor and Tao prove a localization result for microdoubling spaces: One does not
need to consider the supremum over all $r > 0$ when proving weak type bounds, provided the averaging operators
are well behaved. And this is  implied by strong $n$-microdoubling. In the specific case
of $\mathbb{R}^d$ with Lebesgue measure, their localization result entails that it is
enough to consider radii $r$ satisfying $1 \le r\le d$. It is clear that localized maximal operators with
$c\le r \le (1 + 1/d) c$, are  bounded by the averaging operator with radius $r = c$ times
the microdoubling constant. Since $(1 + 1/d)^{d \log d} \approx d$, it follows that we need
roughly $d \log d$ steps to go from 1 to to $d$ by using $c_0 = 1$, $c_1 = (1 + 1/d)$,
$c_2 =  (1 + 1/d)^2$, etc. Thus the maximal operator $M$ with $1 \le r\le d$
is controlled by the sum of $O( d\log d)$ maximal operators with $c_i\le r \le (1 + 1/d) c_i$,
which yields the result by Stein and Str\"{o}mberg mentioned above. Localization is proved by
approximating in a certain sense metric spaces by ultrametric spaces via ``random partitioning methods";
certain modified Doob's maximal inequalities for sublinear operators are proved and applied in their
arguments. A second proof of the $O( d\log d)$ bound is given via the ``Random Vitali Covering Lemma" of
E. Lindenstrauss.

Another setting where it is natural to explore these issues is that of $d$-dimensional Riemannian 
or sub-Riemannian manifolds, or
spaces not as general as metric measure spaces.
In \cite{Li},  Hong-Quan  Li extends to the Heisenberg groups the $O( d)$ estimate of Stein and Str\"{o}mberg 
for euclidean balls on
$\mathbb{R}^d$, by semigroup methods. And in
\cite{LiLo}, Li and Lohou\'e give an $O( d\log d)$ upper bound for the weak type (1,1) inequalities,
when working with the Riemannian volume in hyperbolic
spaces. This is quite remarkable, as the volume of balls in hyperbolic spaces grows exponentially, so
no doubling or microdoubling condition is satisfied (in fact, no doubling measure can be
defined in the hyperbolic spaces). Again the result is obtained by semigroup methods.
In a recent preprint (personal communication) Hong-Quan  Li obtains $L^p$ bounds independent
of the dimension ($p > 1$) for the centered maximal operator in hyperbolic spaces (once
more by semigroup methods). 

Curiously, the analogous question for area on the $d$-dimensional
sphere appears not to have been answered. Of course, one would expect the same result to hold,
that is, the existence of $L^p$ bounds ($p > 1$) independent
of the dimension,  for the centered maximal operator defined by geodesic balls (spherical caps).

A different line of research explores what happens in $\mathbb{R}^d$ under measures that may be different from
Lebesgue measure, restricted to some special class of functions (something which of course, simplifies arguments). 
From now on we always refer to the centered maximal function defined by euclidean balls. It
is shown in \cite[Theorem 3]{MeSo} that considering only radial functions (with Lebesgue measure) leads  to $c_{1,d}\le 4$  in all dimensions, and the same happens if Lebesgue measure is replaced by
a radial, radially increasing measure, cf. \cite[Theorem 2.1]{In}. Besides, for Lebesgue measure and 
radial {\em decreasing} functions, 
it is shown in \cite[Theorem 2.7]{AlPe4} that the sharp constant is $c_{1,d}=1$.

If instead of radial, radially increasing measures one considers radial, radially {\em decreasing} measures,
the situation changes radically. Typically, one has exponential increase in the dimension for $c_{1,d}$,
and some times even for the strong type constants $C_{p, d}$. 
Furthermore it is enough to consider characteristic functions of balls centered at zero 
(hence, radial and decreasing) to prove
exponential increase.
The weak type $(1,1)$ case  for integrable radial
densities defined via bounded decreasing functions was studied in \cite{A1}.
It was shown there that the best constants $c_{1,d}$
satisfy $c_{1,d} \ge \Theta\left(1\right)\left(
2/\sqrt 3\right)^{d/6}$, in strong
contrast with the linear $O(d)$ upper bounds known for Lebesgue measure.
Exponential increase was also shown 
for the same measures and small values of $p > 1$ in \cite{Cri}; shortly after
(and independently) these results were improved in \cite{AlPe5}, as they applied to larger
exponents $p$ and to a wider class of measures. It was also shown in \cite{AlPe5} that
exponential increase could occur for arbitrarily large values of $p$ and suitably chosen
doubling measures. Together with the results for hyperbolic spaces  mentioned before, this
shows that the doubling condition is neither necessary nor sufficient to have ``good
bounds" for maximal inequalities in terms of the dimension. Finally, it is proven in
\cite{CriSjo} that for the standard gaussian measure in $\mathbb{R}^d$, one has 
exponential increase in the constants {\em for all} $p \in (1,\infty)$. So from this viewpoint,
the most important measures in $\mathbb{R}^d$, Lebesgue and Gaussian, behave in a completely
opposite manner.

In the next section we consider the following question about the maximal operator acting on radial
functions: As we have seen, uniform bounds hold for radial non-decreasing measures, and
we have exponential increase for several classes of radial decreasing measures. So it is natural to ask
whether Lebesgue measure is the borderline case which separates uniform from non-uniform behavior in
the constants. We shall show in the next section that the answer to this question is negative:
For the the radial 
decreasing measures $\mu_d$ on $\mathbb{R}^d$, defined by $d\mu_d(y)=\frac{dy}{\|y\|_2^\alpha}$, $\alpha > 0$, and the maximal operator acting on radial integrable functions,
the constants $c_{1,d}$ are bounded uniformly in $d$; of course, the bounds we find increase with
$\alpha$, as was to be expected. In fact, if the exponents $\alpha_d$ are allowed
to increase to infinity with $d$, then so do the constants $c_{1,d}$.

\section{Uniform bounds for some radial  measures and radial functions}

Recall that $\|x\|_2 :=\left( x_1^2+ x_2^2+\dots + x_d^2\right)^{1/2}$.
A function $f:\mathbb{R}^d\to \mathbb{R}$ is radial if there is a second  
function $f_0:(0, \infty) \to \mathbb{R}$
such that 
\begin{equation}\label{defradfunct}
f(x) = f_0(\|x\|_2) 
\end{equation}
on $\mathbb{R}^d\setminus \{0\},$ i.e., 
$f(x)$ depends only on the distance from $x$ to the origin, and
not on $x$ itself (no restriction is placed on $f(0)$). 
Thus,  $f$ is rotation invariant. Since $f$ depends only on one parameter (the
distance to the origin) it is not surprising that uniform bounds can be found 
(at least for some measures) by reduction to the $1$-dimensional case.
All functions considered in this section are radial. Next, radial measures are defined as follows.
Fix $d\in \mathbb{N}\setminus \{0\}$, and
let $\mu_0: (0,\infty ) \to [0,\infty )$ be a (possibly unbounded)
 function, not zero almost everywhere, such that $\mu_0(t) t^{d-1} \in 
 L^1_{\operatorname{loc}}[(0,\infty), dt]$. Then the function $\mu_0$
defines a   rotationally invariant measure
$\mu$  on $\mathbb{R}^d$ via
\begin{equation}\label{defrad}
\mu (A) := \int_A  \mu_0(\|y\|_2)  d\lambda^{d} (y).
\end{equation}
Here $\mu_0$ is allowed to
depend on $d$, and the local integrability of $\mu_0(t) t^{d-1}$ is assumed for each  {\em fixed}
$d$. 
 Furthermore, $\mu$ may fail to be locally finite, even if
$\mu_0(t) t^{d-1} \in 
 L^1_{\operatorname{loc}}[(0,\infty), dt]$. This happens, for instance, if $d = 1$ and 
 $\mu_0(t) = t^{-1}$: In this case $\mu (-h, h) = \infty$ for every $h > 0$. For convenience,
 we assume in this section that maximal operators are defined using closed balls, which
 we denote also by $B(x,r)$, to keep the notation simple.

We shall show next that uniform weak type (1,1) bounds hold for the radial measures with densities 
given  $d\mu(y)=\frac{dy}{\|y\|_2^\alpha}$, where  $\alpha$ is a fixed constant, 
independent of the dimension. However, as soon as we allow the exponents to grow to infinity with
the dimension, this result fails. So the measures $d\mu(y)=\frac{dy}{\|y\|_2^\alpha}$ represent
the borderline case between uniform and non-uniform weak type (1,1) bounds. 
Finally, if the exponents are allowed to grow like $\alpha d$, where $\alpha\in (1/2, 1)$ is
fixed, then there is exponential increase of the constants $C_{p,d}$ {\em for all $p < \infty$}.

\

\begin{theorem}\label{main_theor} 
For $d \ge 1$, let $\mu_{\alpha_d}$ be the measure on $\mathbb{R}^d$ defined by
 $d\mu_{\alpha_d}(x)=\|x\|_2^{-\alpha_d}dx$.  We consider the centered maximal operator
defined by $\mu_{\alpha_d}$ and euclidean balls, acting on radial functions.
 
 1) If the fixed constant $\alpha > 0$ satisfies 
 $1/2<\alpha < 1$ and $\alpha_d := \alpha d$, then for every $p\in [1,\infty)$ there exists a $b = b(p) > 1$
such that $c_{p,d}\ge \Theta (b^d)$. That is, we have exponential increase in the weak type $(p,p)$  bounds for
all $p < \infty$.

2) For  $\alpha_d \le d/2$,    we have $c_{1,d}\ge \Theta ((5^{1/2}/2)^{\alpha_d})$. In particular,  if $\limsup_d \alpha_d = \infty$,
then we always have
$\limsup_d c_{1,d} = \infty$.
 
 3) If  $\sup_d\alpha_d \le \alpha < \infty$,  then there exists a $C = C( \alpha)$
 such that for every $d\ge 1$, $c_{1,d}\le C$. Thus, there are bounds, uniform in the dimension,
 for the weak type (1,1) constants, and hence, by interpolation, for the strong $(p,p)$ constants, 
 whenever
 $1<p<\infty$.
 \end{theorem}
 
 \begin{remark}\label{lezero} If $\alpha_d\le 0$, then we are in the
 case of radial non-decreasing measures, so  $c_{1,d}\le 4$, as we noted above. 
  \end{remark}
  
  \begin{remark} Parts 1) and 3) of the preceding theorem have been independently discovered by
A. Criado in his Ph. D. Thesis, cf. \cite{Crith}. Remarkably, it is also shown there   that Stein's result regarding strong $L^p$ bounds uniform in $d$, for euclidean balls and
  Lebesgue measure, extends to the measures $d\mu_{\alpha}(x)=\|x\|_2^{-\alpha}dx$, $\alpha > 0$
  (without restricting the action of the operator to radial functions, as we do here).
  \end{remark}
  
{\em Proof of part 1)} We follow the same steps as in the proof of
\cite[Theorem 2.8]{CriSjo}, with the appropriate modifications. 
Let $B_r:=B(0,r)$, and denote by $\omega_{d-1}=\sigma_{d-1}(\mathbb{S}^{d-1})$ the area of the
unit sphere $\mathbb{S}^{d-1}$ in $\mathbb{R}^{d}$.  

\begin{lemma} \cite[Lemma 3.1]{CriSjo} 
Let $\mu$ be a rotation-invariant locally finite Borel measure in $\mathbb{R}^d$. 
For all $x\in\mathbb{R}^d$ and all $r,R>0$ such that $\mu (B_r), \mu(B(x,R))>0$, we have
\begin{equation*}
c_{\mu,p}\ge M_\mu \chi_{B_r}(x)\left(\frac{\mu(B_{|x|})}{\mu(B_r)}\right)^{1/p}
\ge \frac{\mu(B(x,R)\cap B_r)}{\mu(B(x,R))}\left(\frac{\mu(B_{|x|})}{\mu(B_r)}\right)^{1/p}.
\end{equation*}
\end{lemma}

Let $\mu_d$ be the Radon measure $d\mu_d(x)=\|x\|_2^{-\alpha d}dx$ in $\mathbb{R}^d$. Assume $1/2<\alpha<1$. We
point out that the
arguments below also work if instead of a constant $\alpha$ we use variables 
 $\beta_d$, provided they belong to a compact subinterval of $(1/2,1)$. 
That is, if $\beta_d$ tends to $1/2$, then the base of exponentiation tends to $1$. And if $\beta_d$ tends to $1$, 
some ``constants" appearing below may explode.

In view of the preceding lemma, it is enough to show that for each fixed $\alpha\in (1/2, 1)$, 
there exist  $r\equiv r(\alpha), R\equiv R(\alpha), c \equiv c(\alpha),  C\equiv C(\alpha)>0$ 
with $r, R <1$,
and $a \equiv a(\alpha)>1$, 
such that
\begin{equation}\label{uno}
\frac{\mu_d(B(e_1,R)\cap B_r)}{\mu_d(B(e_1,R))}\ge \frac{c}{\sqrt{d}},
\end{equation}
and
\begin{equation}\label{dos}
\frac{\mu_d(B_{1})}{\mu_d(B_r)}\ge C a^d.
\end{equation}

Integration in spherical coordinates shows  that for all $\rho >0$,
\begin{equation*}
\mu©(B_\rho)=\frac{\omega_{d-1} }{d(1-\alpha)}\rho^{d(1-\alpha)}.
\end{equation*}
Thus,
\begin{equation*}
\frac{\mu_d(B_{1})}{\mu_d(B_r)}\ge \left(\frac{1}{r}\right)^{(1-\alpha)d},
\end{equation*}
and (\ref{dos}) follows with $C=1$ and $a=(1/r)^{1-\alpha}$.

Next we bound $\mu_d(B(e_1,R))$  from above, by changing  to spherical coordinates:
\begin{equation}\label{por_arriba}
\mu_d(B(e_1,R))= \int_{1-R}^{1+R}|\partial B_s\cap B(e_1,R)|_{d-1}  s^{-\alpha d} ds,
\end{equation}
where $|\cdot|_{d-1} $ denotes the $n-1$ dimensional Hausdorff measure. Call $\beta_s$ the angle determined by the segment that joins the origin with $e_1$ and the one that connects the origin to
any point of intersection of $\partial B_s$ with $\partial  B(e_1,R)$. Then $0\le \beta_s<\pi/2$, since $R<1$.
Thus,
\begin{equation}\label{caps}
|\partial B_s\cap B(e_1,R)|_{d-1} =\int_0^{\beta_s} \omega_{d-2} (s\sin \theta)^{d-2}  s d\theta =\omega_{d-2}  s^{d-1}  \int_0^{\beta_s}(\sin \theta)^{d-2}   d\theta.
\end{equation}
By the cosine law, applied to the triangle $T(1, s, R)$ with side lengths 1, $s$, and $R$, and the angle
$\beta_s$ facing the $R$-side, we have
\begin{equation}\label{coseno}
\cos \beta_s=\frac{1+s^2-R^2}{2s},
\end{equation}
so
\begin{equation}\label{seno}
\sin \beta_s=\left[1-\left(\frac{1+s^2-R^2}{2s}\right)^2\right]^{1/2}.
\end{equation}

Note that the maximum value of $\beta_s$ occurs when the ray starting at 0 is tangent
to $B(e_1,R)$, so the triangle $T(1, s, R)$ has a right angle, and hence $s=\sqrt{1-R^2}$.
Since $\sin \beta_s$ increases with $\beta_s$ and $\cos \beta_s$ decreases, 
from (\ref{coseno}) and (\ref{seno}) we obtain  $\cos \beta_s\ge \sqrt{1-R^2}$ and $\sin \beta_s \le R$.

Using (\ref{caps}) we conclude that 
\begin{equation}\label{equivalencia}
\frac{\omega_{d-2} }{d-1}  (s \sin\beta_s)^{d-1} 
 \le |\partial B_s\cap B(e_1,R)|_{d-1} 
 =\omega_{d-2}  s^{d-1}  \int_0^{\beta_s}(\sin \theta)^{d-2}   d\theta
\end{equation}
\begin{equation}\label{equivalencia2}
\le \frac{ \omega_{d-2}  s^{d-1} }{\sqrt{1-R^2}} \int_0^{\beta_s}\cos \theta (\sin \theta)^{d-2}   d\theta
\le \frac{1}{\sqrt{1-R^2}}\frac{\omega_{d-2} }{d-1}  (s \sin\beta_s)^{d-1} . 
\end{equation}
Define 
\begin{equation}\label{efar}
F_R\left(s\right): =(s \sin\beta_s)^2 s^{-2\alpha}=\frac{1}{4}\left[4s^2-\left(1+s^2-R^2\right)^2\right]s^{-2\alpha}.
\end{equation}
By (\ref{equivalencia}) and (\ref{por_arriba}), 
\begin{equation*}
\mu_d(B(e_1,R))\le \frac{1}{\sqrt{1-R^2}}\frac{\omega_{d-2} }{d-1}  \int_{1-R}^{1+R} (s \sin\beta_s)^{d-1} s^{-\alpha d} ds
\end{equation*}
\begin{equation*}
= \frac{1}{\sqrt{1-R^2}}\frac{\omega_{d-2} }{d-1}  \int_{1-R}^{1+R} (s \sin\beta_s)^{d-1} s^{\alpha(1-d)} \frac{ds}{s^\alpha}
\end{equation*}
\begin{equation*}
=\frac{1}{\sqrt{1-R^2}}\frac{\omega_{d-2} }{d-1}  \int_{1-R}^{1+R} F_R\left(s\right)^{\frac{d-1} {2}}  \frac{ds}{s^\alpha}.
\end{equation*}
Clearly, $F_R(1-R) = F_R(1+R)= 0$. Furthermore, $F_R$ is increasing on 
$[1-R, \sqrt{1-R^2}]$ since it is the product of two increasing functions there
($(\sin\beta_s)^2$ and $s^{2-2\alpha})$).

\vskip .2 cm

{\em Claim} (to be proven later): Choosing $R =\sqrt{1-4(1-\alpha)^2}$, the function 
$F_R$ achieves
its unique maximum on $[1-R, 1 + R]$ at a point $s_0 <1$. 

\vskip .2 cm

Assuming the claim, if we replace
$F_R\left(s\right)$ and $s^{-\alpha}$ in the preceding integral by their maximum values, 
 we obtain
\begin{equation}\label{est1}
\mu_d(B(e_1,R))\le\frac{2R}{(1-R)^\alpha\sqrt{1-R^2}} \frac{\omega_{d-2} }{d-1}  F_R\left(s_0\right)^{\frac{d-1} {2}}.
\end{equation}

Next we set $r:=s_0$. To bound   $\mu_n(B(e_1,R)\cap B_r)$ from below, we change to spherical coordinates and use (\ref{equivalencia}): 
\begin{equation*}
\mu_d(B(e_1,R)\cap B_{s_0})=\int_{1-R}^{s_0}|\partial B_s\cap B(e_1,R)|_{d-1}  s^{-\alpha d} ds\ge
\end{equation*}
\begin{equation}\label{por_abajo}
\frac{\omega_{d-2} }{d-1} \int_{1-R}^{s_0} (s \sin\beta_s)^{d-1}   s^{-\alpha d} ds= \frac{\omega_{d-2} }{d-1} 
\int_{1-R}^{s_0}F_R\left(s\right)^{\frac{d-1} {2}}  \frac{ds}{s^\alpha}.
\end{equation}
By Taylor's approximation, 
 for every $s\in [1-R,1+R]$ there exists a $\tau_s$ between $s$ and $s_0$ such that
\begin{equation*}
 F_R(s)=F_R(s_0)+\frac{F''_R(\tau_s)}{2}(s-s_0)^2.
\end{equation*}
 Denote by $M\equiv M(\alpha)$ the maximum value of $|F_R''|$ on $[1-R,1+R]$. 
 We assume that $d>>1$ is so large that 
 $$
 0<\delta:=\sqrt{4F_R(s_0)/M(d-1)} < s_0 - 1 + R
 $$
 (we can do this since neither $R$ nor $F_R$ depend on $d$).
  Then, for all $s\in (s_0-\delta,s_0)$,
\begin{equation*}
F_R(s)\ge F_R(s_0)-\frac{M}{2}\delta^2=F_R(s_0)\left(1-\frac{2}{(d-1)}\right).\end{equation*}
Since $(1 - t)^{1/t}$ increases to $1/e$ as $t\downarrow 0$, for all $d\ge 4$.
\begin{equation*}
F_R(s)^{\frac{d-1} {2}}\ge F_R(s_0)^\frac{d-1} {2}\left(1-\frac{2}{(d-1)}\right)^\frac{d-1} {2}
\ge F_R(s_0)^\frac{d-1} {2}\left(\frac{1}{3}\right)^\frac{3}{2}.
\end{equation*}
Thus, by (\ref{por_abajo})
\begin{equation*}
\mu_d(B(e_1,R)\cap B_{s_0})\ge \frac{\omega_{d-2} }{d-1} \int_{1-R}^{s_0}F_R\left(s\right)^{\frac{d-1} {2}}  \frac{ds}{s^\alpha}
\end{equation*}
\begin{equation*}
\ge 
\frac{\omega_{d-2} }{d-1} \int_{s_0-\delta}^{s_0}F_R\left(s\right)^{\frac{d-1} {2}} 
 \frac{ds}{s^\alpha}
\end{equation*}
\begin{equation*}
\ge \frac{\omega_{d-2} }{d-1} F_R\left(s_0\right)^{\frac{d-1} {2}}\left(\frac{1}{3}\right)^\frac{3}{2}
\int_{s_0-\delta}^{s_0} \frac{ds}{s^\alpha}
\end{equation*}
\begin{equation}\label{est_2}
\ge \frac{\omega_{d-2} }{d-1} F_R\left(s_0\right)^{\frac{d-1} {2}}\left(\frac{1}{3}\right)^\frac{3}{2}
s_0^{-\alpha} \delta.
\end{equation}
Finally, using (\ref{est1}) and (\ref{est_2}), we get 
\begin{equation*}
\frac{\mu_n(B(e_1,R)\cap B_{s_0})}{\mu_n(B(e_1,R))}\ge 
\frac{
\left(\frac{1}{3}\right)^\frac{3}{2}(1-R)^\alpha\sqrt{1-R^2}s_0^{-\alpha} \delta}
{2R}
\ge \frac{c}{\sqrt{d}},
\end{equation*}
where $c = c(\alpha) > 0$ ($c$ depends on $R$, but recall that $R =\sqrt{1-4(1-\alpha)^2}$).

\

{\em Proof of the claim.} For simplicity, we make the change of variables $t = s^{2}$,
and  write 
 \begin{equation}\label{gee}
g(t) := 4 F_R\left(t^{1/2}\right) = \left[4t-\left(1+t-R^2\right)^2\right]t^{-\alpha}.
\end{equation}
Clearly it is enough to show that $g$ has a unique maximum  
$t_0\in \left[(1-R)^2, (1 + R)^2  \right]$ 
such that $t_0 < 1$. It then follows that $F_R$ has a unique maximum 
$s_0\in [1-R, 1 + R]$  with  
 $s_0 = t_0^{1/2} < 1$.
 
Replacing $R^2$ by its value $1-4(1-\alpha)^2$
in (\ref{gee}) and simplifying we obtain
 \begin{equation}\label{gee1}
g(t) =  \left[-16(\alpha - 1)^4 + (-4 + 16 \alpha - 8 \alpha^2) t-t^2\right]t^{-\alpha}.
\end{equation}
To find the local extrema we differentiate and rearrange:
\begin{equation}\label{geeprime}
g^\prime (t) =  \left[16(\alpha - 1)^4 \alpha + (-4 + 20 \alpha - 24 \alpha^2 + 8\alpha^3) t+ (\alpha-2) t^2\right]
/t^{1+\alpha}.
\end{equation}
Note that  the zeroes of $g^\prime$ are the same as the zeroes of its numerator, so by solving  a second degree
equation, we get
$$t_0 = 4(\alpha- \alpha^2) \mbox{ \ \ \ and \ \ \ } t_1 = \frac{4(\alpha - 1)^3}{2 - \alpha}. 
$$
Now at least one root belongs to $\left[(1-R)^2, (1 + R)^2  \right]$,
since $g$ vanishes at the endpoints and it must have a global maximum.
But $t_1 < 0$,
so the only solution in $\left[(1-R)^2, (1 + R)^2  \right]$ is $t_0$,
and thus the global maximum of $g$ occurs there. Furthermore, on $(1/2, 1)$,
 $f(\alpha) := \alpha - \alpha^2 < 1/4$, whence $t_0=t_0(\alpha) < 1$.
 
This finishes the proof of Part 1).
\qed

\vskip .2 cm

 {\em Proof of part 2).}  Assume that $0<\alpha_d\le d/2$. 
It is shown next  that if $d\ge 12$,
  then
\begin{equation*}
c_{1,d}\ge \frac{1}{2e}\left(\frac{5}{4}\right)^\frac{\alpha_d}{2} .
\end{equation*}

The   proof we present below illustrates the discretization technique, 
valid only for $p=1$.
In this particular application, a radial decreasing function is replaced by one Dirac delta at the
origin. Clearly, any lower bound obtained using $\delta_0$ can be
approximated as much as we want, by considering instead 
the function $\chi_{B(0,r)}/\mu_{\alpha_d}(B(0,r))$,
where $0 < r << 1$. In fact, by the 1-homogeneity of the operator, we can just take $\chi_{B(0,r)}$, since
constants  cancel out. We note that the proofs of exponential growth of the weak and strong type constants 
 in the
papers  
\cite{A1},
\cite{AlPe5}, 
\cite{Cri},
\cite{CriSjo}, all use this method of considering 
$\delta_0$ or $\chi_{B(0,r)}$, and then estimating how shifting balls away from the origin reduces
their measure (the differences between these papers lie in the values of $r>0$ selected, the shifted balls chosen, 
and how their sizes are controlled).

We utilize the following special case of \cite[Proposition 2.1]{A1}: 
\begin{equation}\label{bull}
c_{1,d} \ge 
\frac{\mu_{\alpha_d}( B(0, 1))}{\mu_{\alpha_d} (B(e_1, 1))},
\end{equation}
where $e_1$ is the first vector in the
standard basis of $\mathbb{R}^d$ (any vector of length one will do, by rotational invariance). 
This lower bound is obtained by noticing that 
$M_{\mu_{\alpha_d}} \delta_0 (x) = 1/\mu_{\alpha_d} (B(x,\|x\|_2)$ (recall that balls
can be taken to be closed) and that 
$$
B(0, 1) \subset\left \{M_{\mu_{\alpha_d}} \delta_0 
\ge \frac{1}{\mu_{\alpha_d}(B(e_1, 1))}\right \}.
$$

So, all we need to do is to estimate from below the quotient appearing in (\ref{bull}).
Writing $\sigma_{d-1}$ for the $(d-1)$-dimensional Hausdorff measure on $\mathbb{S}^{d-1}$ (the unit sphere in $\mathbb{R}^d$) integration in polar coordinates yields
\begin{equation}\label{bola_centrada}
\mu_{\alpha_d}(B(0,1))=\frac{\sigma_{d-1}(\mathbb{S}^{d-1})}{d-\alpha_d}.
\end{equation}

Next, note that $B(e_1,1)$ can be decomposed in vertical sections as follows:
\begin{equation*}
B(e_1,1)=\{x=(x_1,\ldots,x_d)\in\mathbb{R}^d: \|x-e_1\|_2\le 1\}=
\end{equation*}
\begin{equation*}
 \{x : 0\le x_1\le 2, (x_2,\ldots,x_d)\in\mathbb{R}^{d-1}, x_2^2+\ldots +x_d^2\le 2x_1-x_1^2\}.
\end{equation*}
Thus, by Fubini's theorem,
\begin{equation*}
\mu_{\alpha_d}(B(e_1,1))=\int_{B(e_1,1)} \frac{dx}{\|x\|_2^{\alpha_d}}=
\end{equation*}
\begin{equation*}
=\int_0^2 \left(\int_{\{(x_2,\ldots,x_d)\in\mathbb{R}^{d-1}, x_2^2+\ldots +x_d^2\le 2x_1-x_1^2\}}\frac{1}{(x_1^2+x_2^2+\ldots +x_d^2)^{\alpha_d/2}}dx_2\cdots dx_d\right)dx_1
\end{equation*}
\begin{equation*}
=:\int_0^2 F(x_1)dx_1,
\end{equation*}
where $F(x_1)$ denotes the inner integral.
Using a spherical change of coordinates we get
\begin{equation*}
F(x_1)= \sigma^{d-2}(\mathbb{S}^{d-2})\int_0^{\sqrt{2x_1-x_1^2}}\frac{t^{d-2} dt}{(x_1^2 +t^2)^{\alpha_d/2}}.
\end{equation*}
Thus
\begin{equation*}
\mu_{\alpha_d}(B(e_1,1))= \sigma^{d-2}(\mathbb{S}^{d-2})\int_0^2\left(\int_0^{\sqrt{2x_1-x_1^2}}\frac{t^{d-2} dt}{(x_1^2 +t^2)^{\alpha_d/2}}\right)dx_1.
\end{equation*}
Note that the region of integration in the above expression 
is the upper semicircle centered at $x_1=1$, $t=0$, in the 
$x_1t$-plane.

Hence, by changing to polar coordinates  we obtain
\begin{equation*}
\mu_{\alpha_d}(B(e_1,1))= \sigma^{d-2}(\mathbb{S}^{d-2}) \int_0^{\pi/2}\left(\int_0^{2\cos\theta} \frac{(\rho\sin\theta)^{d-2}\rho}{\rho^{\alpha_d}}d\rho\right)d\theta=
\end{equation*}
\begin{equation*}
=\frac{\sigma^{d-2}(\mathbb{S}^{d-2})}{d-\alpha_d}
\int_0^{\pi/2}(\sin\theta)^{d-2}(2\cos\theta)^{d-\alpha_d}d\theta
\end{equation*}
\begin{equation}\label{bola_descentrada}
=\frac{2^{d-\alpha_d-1}\sigma_{d-2}(\mathbb{S}^{d-2})\beta(\frac{d-\alpha_d+1}{2},\frac{d-1}{2})}{d-\alpha_d}
\end{equation}
By (\ref{bull}), (\ref{bola_centrada}) and (\ref{bola_descentrada}),
\begin{equation}\label{estimaciondegammas}
c_{1,d}\ge \frac{\sigma_{d-1}(\mathbb{S}^{d-1})}{2^{d-\alpha_d-1}\sigma_{d-2}(\mathbb{S}^{d-2})\beta(\frac{d-\alpha_d+1}{2},\frac{d-1}{2})}
=\sqrt{\pi}\frac{\Gamma(\frac{2d-\alpha_d}{2})}{2^{d-\alpha_d-1}\Gamma(\frac{d}{2})
\Gamma(\frac{d-\alpha_d+1}{2})}
\end{equation}
Now we use the Stirling representation of the Gamma function \cite[p.257, 6.1.38]{A}: For every $x>0$, there exists a $\theta\equiv\theta(x)\in [0,1]$ such that
\begin{equation*}
\Gamma(x+1)= \sqrt{2\pi}x^{x+1/2}e^{-x+\theta/(12x)}.
\end{equation*}
Thus, for $d\ge 3$, we have
\begin{equation}\label{cota1gamma}
\Gamma\left(\frac{d}{2}\right)\le e^{1/6}\sqrt{2\pi}\left(\frac{d-2}{2}\right)^{\frac{d-1}{2}}e^{-\frac{d-2}{2}}.
\end{equation}
and
\begin{equation}\label{cota2gamma}
\Gamma\left(\frac{d-\alpha_d+1}{2}\right)\le e^{1/3}\sqrt{2\pi}\left(\frac{d-\alpha_d-1}{2}\right)^{\frac{d-\alpha_d}{2}}e^{-\frac{d-\alpha_d-1}{2}}.
\end{equation}
We also obtain
\begin{equation}\label{cota3gamma}
\Gamma\left(\frac{2d-\alpha_d}{2}\right)\ge \sqrt{2\pi}\left(\frac{2d-\alpha_d-2}{2}\right)^\frac{2d-\alpha_d-1}{2}e^{-\frac{2d-\alpha_d-2}{2}}.
\end{equation}
Using (\ref{estimaciondegammas}), (\ref{cota1gamma}), (\ref{cota2gamma}) and (\ref{cota3gamma}), we get
\begin{equation*}
c_{1,d}\ge \frac{\sqrt{2}}{e}\frac{\left(2d-\alpha_d-2\right)^\frac{2d-\alpha_d-1}{2}}
{2^{d-\alpha_d}\left(d-2\right)^{\frac{d-1}{2}} \left(d-\alpha_d-1\right)^{\frac{d-\alpha_d}{2}}}.
\end{equation*}
Finally, since $d\ge 12$ and $\alpha_d\le d/2$,  
$$
4[4 \left(d-\alpha_d-1\right)] \ge 5(2d-\alpha_d-2),
$$
and 
$$
\left(2d-\alpha_d-2\right)^2\ge 4\left(d-2\right)\left(d-\alpha_d-1\right).
$$
Thus
\begin{equation*}
c_{1,d}\ge \left(\frac{1}{2e}\right)\frac{\left(2d-\alpha_d-2\right)^\frac{2d-\alpha_d}{2}}
{2^{d-\alpha_d}\left(d-2\right)^{\frac{d}{2}} \left(d-\alpha_d-1\right)^{\frac{d-\alpha_d}{2}}}
\end{equation*}
\begin{equation*}
= \frac{1}{2e}
\left(\frac{\left(2d-\alpha_d-2\right)^2}{4\left(d-2\right)\left(d-\alpha_d-1\right)}\right)^{d/2}
\left(\frac{4 \left(d-\alpha_d-1\right)}{2d-\alpha_d-2}\right)^\frac{\alpha_d}{2}
\ge \frac{1}{2e}\left(\frac{5}{4}\right)^\frac{\alpha_d}{2} .
\end{equation*}
 \qed

  \
  
  Regarding part 3), the rest of this paper presents its proof in detail. 
  Since the upper bounds  we obtain increase with the constant $\alpha$
  (cf. Corollary \ref{coro} below)
  the case where $\alpha_d = \alpha$ for all $d\ge 1$ entails the case $\alpha_d \le \alpha$,
  so from now on we suppose that $\alpha_d = \alpha$ for all $d$.
  
  Note  that
 if $d\le \alpha$, then $\mu_{\alpha_d}$ is not locally finite at the origin, so we want to
 allow this possibility in the definitions. Since below $2 \alpha$ there are only finitely
 many dimensions $1,\dots, [2\alpha]$,  to obtain a uniform bound, it is enough to prove that it exists for 
 $d \ge 2 \alpha$, and then 
  take the largest of these (at most) $1 + [2\alpha] $  constants. The case  $d \ge 2 \alpha$ is considered
  in Corollary \ref{coro} at the end of this paper. This corollary   follows from  Theorem \ref{main_teo},
which is obtained by isolating
the property that makes the proofs of \cite[Theorem 2.1]{In}  and \cite[Theorem 3]{MeSo}  work: 
To
each ball, the argument associates a second ball with the same radius, and center nearer to the
origin (perhaps the origin itself). It is enough to assume that this second ball is not
much larger than the first. 

The following (uniform in the dimension) weak type (1,1) inequality was proven in  
\cite[Theorem 2.1]{In} (cf.  \cite[Theorem 3]{MeSo} for Lebesgue measure):
If $M_\mu$ is the maximal operator associated to centered euclidean balls in $\mathbb{R}^d$ with a radial non-decreasing
measure 
$\mu$, 
then for every $t>0$ and every 
radial $f\in L^1$, 
\begin{equation}\label{MS}
t \mu\{Mf > t \} \le 4 \|f\|_1.
\end{equation}
Even though this proof has already appeared in print (save for some trivial modifications)  we include it here because of its didactic value, 
as it illustrates two basic techniques in the
subject: 1) Control a maximal operator in terms of another operator with known bounds.
2) Instead of integrating over a ball, integrate
over a larger (but not much larger) set (perhaps, just a larger ball). 

Regarding 1), the controlling operator will be the one-dimensional, uncentered Hardy-Littlewood
maximal operator. Its boundedness (cf. the next result)
 hinges
upon the fact that from a finite collection of intervals,   two
disjoint subcollections can be extracted, so that their union is the same as the union of the original
collection (as far as we know, this was published first in \cite{Ra}; it seems to
have been rediscovered, as some authors attribute it to Young).
To see why this is true, 
first throw away unnecessary intervals, those contained in the union of the
others, so
no point belongs to three of them; then label the intervals in increasing order,
say, of the left endpoints, and  notice that the subcollections of intervals
with even and with odd indices are disjoint. As a consequence, one immediately 
obtains the next theorem, cf., for instance, \cite{CaFa} (which makes the unnecessary assumption that compact
sets have finite measure) or \cite{A3}. The result is
valid for completely arbitrary Borel measures (countably additive, non-negative and
not identically 0).

Given a Borel measure  $\nu$, we always assume that it has been completed, i.e., that it
has been extended to the $\sigma$-algebra generated by the Borel sets and the sets of
$\nu$-outer measure zero; we also use $\nu$ to denote this extension. 
While the next result is usually stated for the real line, the
same proof works for subintervals. Alternatively, one can consider $\nu$ defined on a subinterval
$I\subset \mathbb{R}$, and extend it to $\mathbb{R}$ by setting $\nu(I^c) = 0$, thus reducing the case of an
arbitrary interval $I$  to the case $I = \mathbb{R}$. In fact, we will only need the particular interval
 $I = (0,\infty)$.

\begin{theorem}\label{1d} Let $\mu$ be a Borel measure on an interval 
$I\subset \mathbb{R}$,  let
$f\in L^1(\mu)$, and let $M_\mu^u$ be the uncentered maximal operator. 
Then for every $\lambda > 0$, 
 \begin{equation}\label{hipotesis}
   \lambda \mu \{M_\mu^u f > \lambda\} \le 2 \|f\|_1.
 \end{equation}
\end{theorem}

\begin{theorem}\label{main_teo}
Let $\mu$ a radial measure on $\mathbb{R}^d$. Suppose there exists a $C>0$ such that
for all $x\in\mathbb{R}^d$ and all $r$ with $0<r\le 1$, we have
\begin{equation}
\label{hipo_teo}
\mu(B(x\sqrt{1-r^2},r\|x\|_2))\le C \mu(B(x,r\|x\|_2)).
\end{equation}
Then, for every radial function $f\in L^1(\mathbb{R}^d,d\mu)$ and every $\lambda>0$,
\begin{equation}
\lambda \mu\{M_\mu f>\lambda\}\le2(C+1)\|f\|_{L^1(\mathbb{R}^d,d\mu)}.
\end{equation}
\end{theorem}

\begin{remark}
Obviously, all radial non-decreasing measures (including the Lebesgue $d$-dimensional measure) satisfy
 condition (\ref{hipo_teo}) with $C=1$, since the size of balls does not increase when they are shifted
 towards the origin. Note also that when $r = 1$, condition (\ref{hipo_teo}) simply says that 
 $\mu(B(0, \|x\|_2))\le C \mu(B(x, \|x\|_2))$.
\end{remark}

{\em Proof of Theorem \ref{main_teo}.}  Since 
$\mu$ is radial, the local integrability of $\mu_0(t) t^{d-1}$ on $(0, \infty)$
together with condition (\ref{hipo_teo}) entail that all balls have finite measure, so we can assume that balls 
$B(y,s)$ are closed.
Let $r > 0$.
The idea is to show that for every $x\in\mathbb{R}^d\setminus \{0\}$ and every ball $B = B(x,r\|x\|_2)$,
the averages $\frac{1}{\mu(B)}\int_B fd\mu$ are  pointwise bounded by the one-dimensional uncentered
maximal function evaluated at $\|x\|_2$, times a certain constant (since the set $\{0\}$ has measure zero, we
can just forget about it; alternatively,  we note that the set $D$ defined below equals $B$ when $x=0$, 
and then the result is immediate).

We prove the pointwise bound by passing to spherical coordinates. Let 
$v$ be a unit vector
such that the ray $\{t(v):t \ge 0\}$ intersects $B$; in what follows, rays will be denoted just
by $t(v)$. If the segment $I$ resulting from this intersection 
contains $\|x\|_2 v$, then we can use the uncentered operator evaluated at $\|x\|_2 v$, and
there is no need to do anything. However,
it may happen that $I$ does not 
contain $\|x\|_2 v$. If so, we enlarge $I$ up to   $\|x\|_2 v$,
and define $D$ to be the union with $B$ of all these enlarged segments. 
Now 
if $r > 1$, then $D=  B(0 ,\|x\|_2))\cup B(x,r\|x\|_2)$,
whence 
$$
\mu (D)  \le \mu B(0,\|x\|_2) + \mu  B(x,r\|x\|_2) 
$$
$$
\le C \mu  B(x,\|x\|_2)  +  \mu  B(x,r\|x\|_2) 
 \le (C + 1) \mu  B(x,r\|x\|_2).
$$
We show next  that
if $r \le 1$, then
$D \subset B(x\sqrt{1-r^2},r\|x\|_2)\cup B(x,r\|x\|_2)$, so 
$$
\mu (D)  \le \mu B(x\sqrt{1-r^2},r\|x\|_2)) + \mu  B(x,r\|x\|_2) \le (C + 1) \mu  B(x,r\|x\|_2)
$$
(thus, in both cases the measure of $D$ is comparable 
to the measure of $B$). 

For each unit vector $v$ 
such that the ray $t(v)$ intersects $D$, let the segment $[a(v), b(v)]$ denote this intersection.
That is, $a(v)$ is the point of entry (of first intersection) of the ray $t(v)$ in $B$ 
(or equivalently, in $D$), and $b(v)$, the point of exit of $D$, i.e., either $b(v)$ is the point
of exit of the ball, or $b(v) = \|x\|_2 v$, whichever is larger. 

Suppose next that the angle between two given  unit vectors $u$, $w$, is acute ($\le \pi/2$), and let $s >0$.
Let $R$ be the length of the segment joining $su$ with its perpendicular projection over the
segment $[0, sw]$. Then $R$ is also  the length of the segment joining $sw$, with its perpendicular projection over 
$[0, su]$. This observation proves that $D\subset B(x\sqrt{1-r^2},r\|x\|_2)\cup B(x,r\|x\|_2)$, 
as follows.
Consider the vector $x$, and let $v$ be any unit vector 
such that the ray $t(v)$ is tangent to $B = B(x,r\|x\|_2)$. Call this point of tangency $t_0(v)$, and note that
the segment from $x$ to $t_0(v)$ is perpendicular to the ray $t(v)$. We use $T$ to denote the set of all unit 
vectors with rays tangent to $B$, and $S$ the set of all unit 
vectors with rays intersecting $B$ (in particular, $T\subset S$).

The
observation above, with $x = su$, $\|x\|_2 v =  sw$, and $R = r\|x\|_2$, shows that
the points in $D\setminus B$ farthest away from the ray $t x$, i.e., the points of
the form $\|x\|_2 v$, are at distance  
$r\|x\|_2$ from their perpendicular projections  over $[0, x]$. These perpendicular
projections equal $x\sqrt{1-r^2}$ by the Pythagorean Theorem, so all the points $\|x\|_2 v$, $v\in T$,
belong to $B(x\sqrt{1-r^2},r\|x\|_2)$. The points $t_0(v)$ are in $B$, so they are also in $B(x\sqrt{1-r^2},r\|x\|_2)$, 
since the latter ball is just $B$ displaced towards the origin. By convexity, 
 the segments $[t_0(v), \|x\|_2 v]$ are fully
contained in $B(x\sqrt{1-r^2},r\|x\|_2)$. This proves that $D\setminus B \subset B(x\sqrt{1-r^2},r\|x\|_2)$,
as desired.

Now, in order to obtain the pointwise bound
\begin{equation}\label{control}
M_\mu f(x)= M_\mu f_0(\|x\|_2) \le (C+1)M^u_{\gamma_0}f_0(\|x\|_2),
\end{equation}
 all we have to do is to  average $f$ over $D$ instead of $B$. Writing the integral in polar (spherical) coordinates,
the averages of  a function over any segment are always controlled by the uncentered one-dimensional 
maximal operator,
evaluated at {\em any} point of the segment. Since all segments in $D$ contain a point of the form $\|x\|_2 v$ (where
$\|v\|_2  =1$) and since
both the measure and the function are radial,  by evaluating the one-dimensional maximal operator always at
the same point $\|x\|_2$, we are actually averaging a constant function, so we get the same value back. We present
the details next.

Recalling the notation from (\ref{defradfunct}) and (\ref{defrad}), let us define the measure $\gamma_0$  on
$(0,\infty)$  via   
 $d\gamma_0(t):=\mu_0(t)t^{d-1}dt$, so given any subinterval $I\subset (0,\infty)$, 
 \begin{equation}\label{interv}
\gamma_0(I) = \int_{I}\mu_0(t)t^{d-1}dt.
\end{equation}
Writing $\sigma$ for area on the unit sphere, and integrating in spherical coordinates, we get
\begin{equation}\label{one}
\frac{1}{\mu(B(x,r\|x\|_2))}\int_{B(x,r\|x\|_2)} |f(y)|d\mu(y)
= 
\frac{\mu(D)}{\mu(B)} \frac{1}{\mu(D)}\int_B |f(y)|d\mu(y)
\le 
\frac{C+1}{\mu(D)}\int_D |f(y)|d\mu(y)
\end{equation}
\begin{equation}\label{tres}
=
\frac{C+1}{\mu(D)}\int_D |f_0(\|y\|_2)|\mu_0(\|y\|_2)dy
= \frac{C+1}{\mu(D)}\int_{S}\left(\int_{a(v)}^{b(v)} |f_0(t)|\mu_0(t)t^{d-1}dt\right)d\sigma(v)
\end{equation}
\begin{equation}\label{cuatro}
\le 
 \frac{C+1}{\mu(D)}\int_{S} \gamma_0((a(v), b(v) ))M^u_{\gamma_0}f_0(\|x\|_2)d\sigma(v) 
 = (C+1) M^u_{\gamma_0}f_0(\|x\|_2).
 \end{equation}
Taking the supremum over $r > 0$ in (\ref{one}), we obtain (\ref{control}). Finally,
we express the level sets of $M_\mu f$ in spherical coordinates, and apply Theorem \ref{1d}:
\begin{equation*}
\mu\{x\in \mathbb{R}^d: M_\mu f(x)>\lambda\}
\le \mu\left\{x\in \mathbb{R}^d: M^u_{\gamma_0} f_0(\|x\|_2)>\frac{\lambda}{C+1}\right\}
\end{equation*}
\begin{equation*}
= \int_{\mathbb{S}^{d-1}}\left(\int_{\left\{M^u_{\gamma_0} f_0>\frac{\lambda}{C+1}\right\}}\mu_0(t)t^{d-1}dt\right)d\sigma(\omega)
=
\int_{\mathbb{S}^{d-1}}
\gamma_0\left\{M^u_{\gamma_0} f_0>\frac{\lambda}{C+1}\right\}d\sigma(\omega)
\end{equation*}
\begin{equation*}
 \le \frac{2(C+1)}{\lambda}\int_{\mathbb{S}^{d-1}}\left(\int_{(0,\infty)} |f_0(t)|d\gamma_0(t)\right)d\sigma(\omega)= \frac{2(C+1)}{\lambda}\int_{\mathbb{R}^d}|f|d\mu. 
\end{equation*}
\qed

To bound $\mu(B(x,r\|x\|_2))$ from below in the next result, in expressions (\ref{cor_est_sup0}) and 
(\ref{cor_est_sup}) below, 
it is enough to replace the density by its lowest value on $B(x,r\|x\|_2)$, that is, by  $(\|x\|_2 (1+r))^{-\alpha}$.
We use $\left(\|x\|_2\sqrt{1+r^2}\right)^{-\alpha}$ instead, noting that the density is larger than this constant on
at least half the ball. The estimates are not very different, but the second choice gives better
constants for high values of $\alpha$.

\begin{corollary}\label{coro}
Fix  $\alpha>0$ and set $d\mu(y)=\frac{dy}{\|y\|_2^\alpha}$ on $\mathbb{R}^d$, for $d \ge 1$.
If $d\ge 2\alpha$ and $f\in L^1(\mathbb{R}^d,d\mu)$ is radial, then for every $\lambda>0$,
\begin{equation*}
\mu\{M_\mu f>\lambda\}\le \frac{2(4\cdot 6^{\alpha/2}+1)}{\lambda}\|f\|_{L^1(\mathbb{R}^d,d\mu)}.
\end{equation*} 
\end{corollary}

\begin{proof}
We  show that (\ref{hipo_teo}) holds with $C=4\cdot 6^{\alpha/2}$.
Because $\mu$ is radial decreasing, it is clear that the measure of balls increases when they are
shifted towards the origin, since the density is always larger on all points of the shifted ball that are not
contained in the intersection (of the two balls) than on the points of the original ball not in the intersection.
Thus
\begin{equation}\label{cor_est1}
\mu(B(x\sqrt{1-r^2},r\|x\|_2))\le \mu(B(0,r\|x\|_2))= \frac{d}{d-\alpha} (\|x\|_2r)^{d-\alpha} v_d \le 2 (\|x\|_2r)^{d-\alpha} v_d , 
\end{equation}
where $v_d$ denotes the Lebesgue $d$-dimensional measure of the unit ball.

On the other hand,
\begin{equation}\label{cor_est_sup0}
\mu(B(x,r\|x\|_2))=\int_{B(x,r\|x\|_2)}\frac{dy}{\|y\|_2^\alpha}\ge \int_{B(x,r\|x\|_2)\cap \{y:\|y\|_2\le \|x\|_2\sqrt{1+r^2}\}}
\frac{dy}{\|y\|_2^\alpha}\ge
\end{equation}
\begin{equation}\label{cor_est_sup}
 \frac{1}{\|x\|_2^\alpha(1+r^2)^{\alpha/2}} \lambda^d(B(x,r\|x\|_2)\cap \{y:\|y\|_2\le \|x\|_2\sqrt{1+r^2}\})
\end{equation} 
\begin{equation}\label{cor_est_sup2}
\ge \frac{ \lambda^d(B(x,r\|x\|_2))}{2\|x\|_2^\alpha(1+r^2)^{\alpha/2}}=\frac{(r\|x\|_2)^d \  v_d}{2\|x\|_2^\alpha(1+r^2)^{\alpha/2}}
\end{equation} 

If $1/\sqrt{5}\le r\le 1$, it follows from (\ref{cor_est_sup0}-\ref{cor_est_sup2}) and (\ref{cor_est1}) that in order
to obtain (\ref{hipo_teo}), it is enough to find a $C^{\prime}  > 0$ such that
\begin{equation*}
 2 (\|x\|_2r)^{d-\alpha} \  v_d \le C^{\prime}  \frac{(r\|x\|_2)^d \  v_d}{2\|x\|_2^\alpha(1+r^2)^{\alpha/2}}.
\end{equation*}
Simplifying, we see that  $C^{\prime} = 4\cdot 6^{\alpha/2}$ suffices.

Suppose next that $0<r\le 1/\sqrt{5}$. Then  
\begin{equation*}
\mu(B(x\sqrt{1-r^2},r\|x\|_2))=\int_{B(x\sqrt{1-r^2},r\|x\|_2)}\frac{dy}{\|y\|_2^\alpha}\le \frac{1}{\|x\|_2^\alpha (\sqrt{1-r^2}-r)^\alpha} \int_{B(x\sqrt{1-r^2},r\|x\|_2)}dy=
\end{equation*}
\begin{equation}\label{cor_est2}
 =\frac{(r\|x\|_2)^d \ v_d }{\|x\|_2^\alpha (\sqrt{1-r^2}-r)^\alpha} 
\end{equation}
Arguing as in the previous case, we see that it is enough to find a $C^{\prime\prime}  > 0$ such that
\begin{equation*}
\frac{(r\|x\|_2)^d \ v_d}{\|x\|_2^\alpha (\sqrt{1-r^2}-r)^\alpha} \le C^{\prime\prime}   \frac{(r\|x\|_2)^d \ v_d}{2\|x\|_2^\alpha(1+r^2)^{\alpha/2}}.
\end{equation*}
Simplifying, we see that we can take  $C^{\prime\prime} =  2 \cdot 6^{\alpha/2}$.
Since $C^{\prime} \ge C^{\prime\prime}$, (\ref{hipo_teo}) follows with $C= C^{\prime}$. 
\end{proof}

And with the proof of Corollary \ref{coro}, the proof or Theorem \ref{main_theor}, Part 3, is also finished.


\begin{thebibliography}{WWW}

\bibitem[A]{A} M. Abramowitz and I. A. Stegun (Eds), 
{\em Handbook of Mathematical Functions with Formulas, Graphs, and Mathematical Tables, National Bureau of Standards.} Applied Mathematics
Series 55, 9th printing, Washington, 1970.


\bibitem[A1]{A1} Aldaz, J. M.
{\em Dimension dependency of the weak type
$(1,1)$ bounds for maximal functions associated to finite
radial measures.}  Bull. Lond. Math. Soc. 39 (2007) 203--208.
Available at the Math. ArXiv.

\bibitem[A2]{A2} Aldaz, J. M. {\em The  weak type $(1,1)$ bounds for the maximal function associated to cubes
grow to infinity with the dimension.} 
 Ann. of Math. (2) 173 (2011), no. 2, 1013--1023. Available at the Math. ArXiv.

\bibitem[A3]{A3} Aldaz, J. M. {\em A general covering lemma for the real line.} 
Real Anal. Exchange 17 (1991/92), no. 1, 394--398.

\bibitem[ACP]{ACP} Aldaz, J. M.; Colzani, L.; P\'erez L\'azaro, J. {\em Optimal bounds on the modulus of continuity of 
the uncentered Hardy-Littlewood maximal function.}
 J. Geom. Anal. 22 (2012), no. 1, 132--167.
 
 
\bibitem[AlPe1]{AlPe1} Aldaz, J.M., P\'erez L\'azaro, J. {\em Functions of
bounded variation, the derivative of the one dimensional maximal
function, and applications to inequalities,} Trans. Amer. Math. Soc.
359 (5) (2007), 2443--2461. Available at the Math.
ArXiv.



\bibitem[AlPe2]{AlPe2} Aldaz, J. M.; P\'erez L\'azaro, J. {\em Boundedness and unboundedness
results for some maximal operators
on functions of bounded variation.}   J. Math. An. Appl.  Volume
337, Issue 1, (2008)  130--143. Available at the Math.
ArXiv.

\bibitem[AlPe3]{AlPe3} Aldaz, J. M.; P\'erez L\'azaro, J. {\em Regularity of the Hardy-Littlewood maximal operator on
 block decreasing functions.}  Studia Math. 194 (3) (2009) 253--277. Available at the Math.
ArXiv.

\bibitem[AlPe4]{AlPe4} Aldaz, J. M.; P\'erez L\'azaro, J. {\em The best constant for the
centered maximal operator  on  radial functions.}  Math. Inequal. Appl. 14 (2011), no. 1, 173--179; available at the Math.
ArXiv.

\bibitem[AlPe5]{AlPe5} Aldaz, J. M.; P\'erez L\'azaro, J. {\em Behavior of weak type bounds for high 
dimensional maximal operators defined by certain radial measures.} Positivity 15 (2011), 199--213.
available at the Math.
ArXiv.



\bibitem[Au]{Au} Aubrun, G. {\em Maximal inequality for high-dimensional cubes.} Confluentes Mathematici,
Volume 1, Issue 2, (2009) pp. 169--179,
DOI No: 10.1142/S1793744209000067. Available at the Math. ArXiv.

\bibitem[Bou0]{Bou0} Bourgain, Jean
{\em Estimations de certaines fonctions maximales.} 
C. R. Acad. Sci. Paris Sér. I Math. 301 (1985), no. 10, 499–-502.

\bibitem[Bou1]{Bou1} Bourgain, J. {\em On high-dimensional maximal functions associated to
convex bodies.} Amer. J. Math. 108 (1986), no. 6, 1467--1476.

\bibitem[Bou2]{Bou2} Bourgain, J. {\em On the $L\sp p$-bounds for maximal functions
associated to convex bodies in $R\sp n$}. Israel J. Math. 54
(1986), no. 3, 257--265.

\bibitem[Bou3]{Bou3} Bourgain, J. {\em On dimension free maximal inequalities for convex
symmetric bodies in $R\sp n$}. Geometrical aspects of functional
analysis (1985/86), 168--176, LNM, 1267,
Springer, Berlin, 1987.

\bibitem[Bou4]{Bou4} Bourgain, J. {\em On the Hardy-Littlewood maximal function for the cube.}
arXiv:1212.2661.

\bibitem[Ca]{Ca} Carbery, A. {\em An almost-orthogonality principle with
applications to maximal functions associated to convex bodies.}
Bull. Amer. Math. Soc. (N.S.) 14 (1986), no. 2, 269--273.

\bibitem[CaFa]{CaFa}  Capri, O. N.; Fava, N. A. 
{\em Strong differentiability with respect to product measures.}
 Studia Math. 78 (1984), no. 2, 173--178. 

\bibitem[Cri]{Cri} Criado, A.  {\em On the lack of dimension free estimates in $L^p$ for maximal functions 
associated to radial measures.}
 Proc. Roy. Soc. Edinburgh Sect. A 140 (2010), no. 3, 541--552. Available at the Math.
ArXiv.

\bibitem[CriSjo]{CriSjo} Criado, A., Sj\"ogren, P.  {\em Bounds for maximal functions 
associated to rotational invariant measures in high dimensions.} To appear, J. Geom. Anal.
arXiv:1111.4358.

\bibitem[Crith]{Crith} Criado, A.   {\em Problems of harmonic analysis in high dimensions.} Ph. D. Thesis,
Advisor: F. Soria. Universidad Aut\'onoma de Madrid, 2012.

\bibitem[In]{In} Infante, A. {\em Free-dimensional boundedness of the maximal operator.} Bol. Soc. Mat. Mexicana (3) Vol. 14, 2008.
  
 \bibitem[Li]{Li} Li, Hong-Quan {\em Fonctions maximales centr\'ees de Hardy-Littlewood sur les groupes de Heisenberg.} 
Studia Math. 191 (2009), no. 1. 
 
 \bibitem[LiLo]{LiLo} Li, Hong-Quan, Lohou\'e, No\"el 
{\em Fonction maximale centr\'ees de Hardy-Littlewood sur les espaces hyperboliques.} To appear, Ark. f\"or Mat.

\bibitem[MeSo]{MeSo} Menarguez, T. and Soria, F. {\em On the maximal operator associated to a convex body in $\mathbb{R}^n$.} Collect. Math. 43(1992), no. 3, 243--251. 
 
 \bibitem[Me1]{Me1}  Melas, Antonios D. {\em On the centered Hardy-Littlewood maximal operator.}
 Trans. Amer. Math. Soc. 354 (2002), no. 8, 3263--3273. 
 
\bibitem[Me2]{Me2}   Melas, Antonios D. {\em The best constant for the centered Hardy-Littlewood maximal inequality.} 
Ann. of Math. (2) 157 (2003), no. 2, 647--688. 

\bibitem[Mu]{Mu} M\"{u}ller, D. {\em A geometric bound for maximal functions
associated to convex bodies.} Pacific J. Math. 142 (1990), no. 2,
297--312.

\bibitem[NaTa]{NaTa} Naor, A.; Tao, T. {\em Random martingales and localization of maximal inequalities.} 
 J. Funct. Anal. 259 (2010), no. 3, 731--779.
Available at the Math. ArXiv.

\bibitem[Ra]{Ra} Rado, T. {\em Sur un probleme relatif a un th\'eor\`eme de Vitali.} Fund. Math 11 (1928),
228--229. 


\bibitem[StPNAS]{StPNAS} Stein, Elias M. {\em  Maximal functions. I. Spherical means.}
 Proc. Nat. Acad. Sci. U.S.A. 73 (1976), no. 7, 2174--2175. 
 
\bibitem[St1]{St1} Stein, E. M. {\em The development of square functions in the work of A.
Zygmund. Bull. Amer. Math. Soc.} (N.S.) 7 (1982), no. 2, 359--376.

\bibitem[St2]{St2} Stein, E. M. {\em Three variations on the theme of maximal
 functions.}
 Recent progress in Fourier analysis (El Escorial, 1983), 229--244,
 North-Holland Math. Stud., 111, North-Holland, Amsterdam, 1985.

\bibitem[St3]{St3} Stein, E. M. {\em Harmonic analysis: real-variable methods,
 orthogonality, and oscillatory integrals.} Princeton University Press, Princeton, NJ, 1993.

 \bibitem[StSt]{StSt} Stein, E. M.; Str\"{o}mberg, J. O. {\em Behavior of maximal functions in $R\sp{n}$
 for large $n$.} Ark. Mat. 21 (1983), no. 2, 259--269.

\end{thebibliography}
\end{document}